\documentclass[english]{article}
\usepackage[T1]{fontenc}
\usepackage[latin9]{inputenc}
\usepackage{geometry}
\usepackage{verbatim}
\usepackage{amsmath}
\usepackage{amssymb}
\usepackage{esint}

\makeatletter
\newcommand{\lyxaddress}[1]{
	\par {\raggedright #1
	\vspace{1.4em}
	\noindent\par}
}

\makeatother

\makeatother

\usepackage{babel}

\makeatother

\usepackage{babel}
\begin{document}

\title{An integral transform for quantum amplitudes}

\author{Jack C. Straton}
\maketitle

\lyxaddress{Department of Physics, Portland State University, Portland, OR 97207-0751,
USA}

\lyxaddress{straton@pdx.edu}
\begin{abstract}
The central impediment to reducing multidimensional integrals of transition
amplitudes to analytic form, or at least to a fewer number of integral
dimensions, is the presence of magnitudes of coordinate vector differences
(square roots of polynomials) $|{\bf x}_{1}-{\bf x}_{2}|^{2}=\sqrt{x_{1}^{2}-2x_{1}x_{2}\cos\theta+x_{2}^{2}}$
in disjoint products of functions. Fourier transforms circumvent this
by introducing a three-dimensional momentum integral for each of those
products, followed in many cases by another set of integral transforms
to move all of the resulting denominators into a single quadratic
form in one denominator whose square my be completed. Gaussian transforms
introduce a one-dimensional integral for each such product while squaring
the square roots of coordinate vector differences and moving them
into an exponential. Addition theorems may also be used for this purpose,
and sometimes direct integration is even possible. Each method has
its strengths and weaknesses. An alternative integral transform to
Fourier transforms and Gaussian transforms is derived herein and utilized.
A number of consequent integrals of Macdonald functions, hypergeometric
functions, and Meijer G-functions with complicated arguments is given. 
\end{abstract}
\vspace{2pc}
 \\

\section{Introduction}

The analytical reduction of atomic integrals involving explicit functions
of the inter-electron distances is the central task for evaluating
transition amplitudes. Direct
integration is sometimes possible (see, for instance, \cite{Ley-Koo and Bunge}, among many
others), and at other times Fourier transforms (e.g, \cite{Fromm and Hill},\cite{Remiddi},
and \cite{Harris PRA 55 1820}), Gaussian transforms (e.g., \cite{Kikuchi},\cite{Shavitt and Karplus},
and \cite{Stra89a}), and addition theorems (e.g. \cite{Sack}, \cite{Porras and  King},
and \cite{Weniger_two-range}) are more useful.

The main drawback of integral transforms is that one must introduce
additional integral dimensions in order to remove the initial ones,
and the reduction of those introduced integrals becomes more difficult
the larger the numbers of wave functions transformed. For Fourier
transforms, one must introduce a three-dimensional integral for each
wave function and often additional integrals to combine the resulting
momentum denominators into a single denominator so that one can complete
the square in the momenta to allow the angular integrals to be performed.
\cite{Stra89c}
Gaussian transforms, on the other hand, require just a single one-dimensional
integral for each wave function, and the completion of the square
in the coordinate variables can be done in the resulting exponential.
The author nevertheless finds the former approach useful as a check
on the latter.

Since more researchers are familiar with Fourier transforms, let us
explicate these ideas using Gaussian transforms. Consider the a product
of two Slater-type atomic orbitals, the seed function $\psi_{000}$
from which Slater functions \cite{Chen} and Hylleraas powers \cite{Harris PRA 55 1820}
are derived by differentiation. (Known as the Yukawa \cite{Yukawa}
exchange potential in nuclear physics, this function also appears
in plasma physics, where it is known as the Debye-H\"uckel potential,
arising from screened charges \cite{NayekGhoshal} requiring the replacement
of the Coulomb potential by an effective screened potential.\cite{EckerWeizel,Harris}
Such screening of charges also appears in solid-state physics, where
this function is called the Thomas-Fermi potential. In the atomic
physics of negative ions, the radial wave function is given by the
equivalent Macdonald function $\left[R(r)=\frac{C}{\sqrt{r}}K_{1/2}(\eta r)\right]$.\cite{Smirnov2003} This function also appears in the approximate
ground state wave function \cite{GaravelliOliveira} for a hydrogen
atom interacting with hypothesized non-zero-mass photons.\cite{CaccavanoLeung}
For simplicity, the term \emph{Slater orbital} will be used for this
function herein.)

\begin{equation}
S_{1}^{\eta_{1}0\eta_{12}0}\left(0;0,\mathbf{x}_{2}\right)\equiv S_{1}^{\eta_{1}j_{1}\eta_{12}j_{2}}\left(\mathbf{p}_{1};\mathbf{y}_{1},\mathbf{y}_{2}\right)_{p_{1}\rightarrow0,y_{1},\rightarrow0,y_{2}\rightarrow x_{2},j_{1}\rightarrow0,j_{2}\rightarrow0}=\int d^{3}x_{1}\frac{e^{-\eta_{1}x_{1}}}{x_{1}}\frac{e^{-\eta_{12}x_{12}}}{x_{12}}\quad,\label{eq:SVxVxy}
\end{equation}
 where we use the much more general notation of previous work \cite{Stra89a}
in which the short-hand form for shifted coordinates is $\mathbf{x}_{12}=\mathbf{x}_{1}-\mathbf{x}_{2}$,
$\mathbf{p}_{1}$ is a momentum variable within any plane wave associated with
the (first) integration variable, the $\mathbf{y}_{i}$ are coordinates external
to the integration, and the \emph{j}s are defined in the Gaussian
transform\cite{Stra89a} of the generalized Slater orbital: 
\begin{equation}
\begin{array}[t]{ccc}
V^{\eta j}({\bf R}) & = & R^{j-1}e^{-\eta R}=\left(-1\right)^{j}{\displaystyle \frac{d^{j}}{d\eta^{j}}}{\displaystyle \frac{1}{\sqrt{\pi}}}\int_{0}^{\infty}\,d\rho_{3}{\displaystyle \frac{e^{-R^{2}\rho_{3}}e^{-\eta^{2}/4/\rho_{3}}}{\rho_{3}^{\;1/2}}}\;\;\left[\eta\geqq0,\:R>0\right]\quad.\\
 & = & R^{j-1}e^{-\eta R}=\hspace{-0cm}{\displaystyle \frac{1}{2^{j}\sqrt{\pi}}}\int_{0}^{\infty}\,d\rho_{3}{\displaystyle \frac{e^{-R^{2}\rho_{3}}e^{-\eta^{2}/4/\rho_{3}}}{\rho_{3}^{\;\left(j+1\right)/2}}}H_{j}\left({\displaystyle \frac{\eta}{2\sqrt{\rho_{3}}}}\right)\;\;\left[\forall j\geq0\:\mathrm{if\:}\eta>0,\;j=0\:\mathrm{if\:}\eta=0\right]
\end{array}\label{seventeen}
\end{equation}

Then

\begin{eqnarray}
S_{1}^{\eta_{1}0\eta_{12}0}\left(0;0,x_{2}\right) & = & \int d^{3}x_{1}{\displaystyle \frac{1}{\sqrt{\pi}}}\int_{0}^{\infty}\,d\rho_{1}{\displaystyle \frac{e^{-x_{1}^{2}\rho_{1}}e^{-\eta_{1}^{2}/4/\rho_{1}}}{\rho_{1}^{\;1/2}}}{\displaystyle \frac{1}{\sqrt{\pi}}}\int_{0}^{\infty}\,d\rho_{2}{\displaystyle \frac{e^{-x_{12}^{2}\rho_{2}}e^{-\eta_{12}^{2}/4/\rho_{2}}}{\rho_{2}^{\;1/2}}}\nonumber \\
 & = & {\displaystyle \frac{1}{\pi}}\int d^{3}x'_{1}\int_{0}^{\infty}\,d\rho_{1}{\displaystyle \frac{e^{-\eta_{1}^{2}/4/\rho_{1}}}{\rho_{1}^{\;1/2}}}\int_{0}^{\infty}\,d\rho_{2}{\displaystyle \frac{e^{-\eta_{12}^{2}/4/\rho_{2}}}{\rho_{2}^{\;1/2}}}\nonumber \\
 & \times & exp\left(-\left(\rho_{1}+\rho_{2}\right)x'{}_{1}^{2}-\frac{x_{2}^{2}\rho_{1}\rho_{2}}{\rho_{1}+\rho_{2}}\right)\quad,\label{eq:SVxVxyG}
\end{eqnarray}
 where we have not displayed the steps involved in completing the
square in the quadratic form in the integration variable $x_{1}$,
which allows the spatial integral to be done by changing variables
from $\mathbf{x}_{1}$ to $ \hspace{1cm} \mathbf{x}'_{1}=\mathbf{x}_{1}-\frac{\rho_{2}}{\rho_{1}+\rho_{2}}\mathbf{x}_{2}$
with unit Jacobian,\cite{GR5 p. 382 No. 3.461.2}

\begin{eqnarray}
\int e^{-\left(\rho_{1}+\rho_{2}\right)x'{}_{1}^{2}}d^{3}x'{}_{1}=4\pi\int_{0}^{\infty}e^{-\left(\rho_{1}+\rho_{2}\right)x'{}_{1}^{2}}\;x'{}_{1}^{2}dx'_{1}=\frac{4\pi^{1+1/2}}{2^{2}\left(\rho_{1}+\rho_{2}\right)^{3/2}}\quad\left[\rho_{1}+\rho_{2}>0\right] \quad.
\end{eqnarray}
 What remains is

\begin{eqnarray}
S_{1}^{\eta_{1}0\eta_{12}0}\left(0;0,x_{2}\right) & = & \pi^{1/2}\int_{0}^{\infty}\,d\rho_{1}{\displaystyle \frac{e^{-\eta_{1}^{2}/4/\rho_{1}}}{\rho_{1}^{\;1/2}}}\int_{0}^{\infty}\,d\rho_{2}{\displaystyle \frac{e^{-\eta_{12}^{2}/4/\rho_{2}}}{\rho_{2}^{\;1/2}}}\label{eq:threerphos-1}\\
 & \times & \frac{1}{\left(\rho_{1}+\rho_{2}\right)^{3/2}}\exp\left(-\frac{x_{2}^{2}\rho_{1}\rho_{2}}{\rho_{1}+\rho_{2}}\right)\quad.\nonumber 
\end{eqnarray}

Let

\begin{equation}
\tau_{1}={\displaystyle \frac{\rho_{1}}{\rho_{1}+\rho_{2}}}\quad,\label{eq:tau1}
\end{equation}
 then \cite{GR5 p. 384 No. 3.471.9}

\begin{eqnarray}
S_{1}^{\eta_{1}0\eta_{12}0}\left(0;0,x_{2}\right) & = & \pi^{1/2}\int_{0}^{1}\,d\tau{\displaystyle \frac{1}{\tau{}^{1/2}}}\int_{0}^{\infty}\,d\rho_{2}{\displaystyle \frac{1}{\rho_{2}^{\;1/2+1}}}\nonumber \\
 & \times & \exp\left(-x_{2}^{2}\tau\rho_{2}-\left(\eta_{12}^{2}\tau+\eta_{1}^{2}\left(1-\tau\right)\right)/\tau/4/\rho_{2}\right)\nonumber \\
 & = & \pi^{1/2}\int_{0}^{1}\,d\tau{\displaystyle \frac{1}{\tau{}^{1/2}}}\int_{0}^{\infty}\,d\rho_{2}{\displaystyle \frac{1}{\rho_{2}^{\;1/2+1}}}\nonumber \\
 & \times & \exp\left(-x_{2}^{2}\tau\rho_{2}-\left(\eta_{12}^{2}\tau+\eta_{1}^{2}\left(1-\tau\right)\right)/\tau/4/\rho_{2}\right)\nonumber \\
 & = & \pi^{1/2}\int_{0}^{1}\,d\tau\frac{2\sqrt{\pi}e^{-x_{2}\sqrt{\tau\left(\eta_{12}^{2}-\eta_{1}^{2}\right)+\eta_{1}^{2}}}}{\sqrt{\tau\left(\eta_{12}^{2}-\eta_{1}^{2}\right)+\eta_{1}^{2}}}\quad.\label{eq:tau int}
\end{eqnarray}
 Changing variables to

\begin{equation}
\begin{array}{c}
s=\left[\tau\left(\eta_{12}^{2}-\eta_{1}^{2}\right)+\eta_{1}^{2}\right]^{1/2}\end{array}\label{eq:s}
\end{equation}
 allows one to perform the indefinite integration \cite{GR5 p. 111 No. 2.311}

\begin{eqnarray}
S_{1}^{\eta_{1}0\eta_{12}0}\left(0;0,x_{2}\right) & = & 2\pi\int_{\eta_{1}}^{\eta_{12}}\,{\displaystyle \frac{2sds}{\eta_{12}^{2}-\eta_{1}^{2}}}\frac{e^{-x_{2}s}}{s}\nonumber \\
 & = & 4\pi{\displaystyle \frac{1}{\eta_{12}^{2}-\eta_{1}^{2}}}\left.\frac{e^{-x_{2}s}}{-x_{2}}\right|_{\eta_{1}}^{\eta_{12}}\nonumber \\
 & = & \frac{4\pi\left(e^{-\eta_{12}x_{2}}-e^{-\eta_{1}x_{2}}\right)}{x_{2}\left(\eta_{1}^{2}-\eta_{12}^{2}\right)}\quad.\label{eq:finalresult}
\end{eqnarray}

In an attempt to simplify this process by bypassing the latter two
changes of variable, the author was able to correct an integral tabled
in Prudnikov, Brychkov, and Marichev\cite{PBM p. 567 No. 3.1.3.4}
that should have read 
\begin{eqnarray}
\int_{0}^{\infty}\int_{0}^{\infty}\frac{1}{\sqrt{x+y}} & f\left(\frac{xy}{x+y}\right)e^{-px-qy}dx\,dy=\frac{\sqrt{\pi}\left(\sqrt{p}+\sqrt{q}\right)}{\sqrt{pq}}\int_{0}^{\infty}e^{-\left(\sqrt{p}+\sqrt{q}\right)^{2}t}f(t)\,dt.\label{eq:fixed}
\end{eqnarray}
 and generalize it to a wide class of integrals \cite{Atoms 8 53 (2020}:

\begin{eqnarray}
R_{2}\left(n,\,m,\,\nu,\,a,\,b,\,c,\,h,\,j,\,p,\,q\right) & = & \int_{0}^{\infty}\int_{0}^{\infty}\frac{1}{x^{n/2}y^{m/2}\left(x+y\right)^{\nu/2}}f\left(\frac{xy}{x+y}\right)\nonumber \\
 & \times & e^{-\frac{a}{x}-\frac{b}{y}-c\,xy/(x+y)-h\,y/(x+y)-j/(x+y)-px-qy}dx\,dy.\label{eq:I2gen}
\end{eqnarray}

\section{Seeking a simpler transform}

For several years the author has been fascinated with a little-used
integral transform \cite{GR5 p. 649 No. 4.638.2 GR7 p. 615} 
\begin{equation}
\frac{1}{r_{0}^{s-p_{1}}r_{1}^{p_{1}}}=\frac{1}{\Gamma\left(p_{1}\right)}\int_{0}^{\infty}\frac{\zeta_{1}^{p_{1}-1}}{\left(r_{1}\zeta_{1}+r_{0}\right){}^{s}}dt\zeta_{1}\label{eq:GR5 p. 649 No. 4.638.2 GR7 p. 615}
\end{equation}
 that has a tantalizing one-fewer integrals than the Gaussian transform,
while nevertheless moving the coordinate variables into a single quadratic
form whose square may be completed. Its downside, of course, is that
it does not apply to Slater orbitals. We can nevertheless show its
utility by setting $\eta_{12}=0$ in eq. (\ref{eq:SVxVxy}) and use
it (with $p_{1}=1/2$ and $s=1$) to reduce the integral over the
product of one Slater orbital and one Coulomb potential:

\begin{eqnarray}
S_{1}^{\eta_{1}000}\left(0;0,x_{2}\right) & = & \int d^{3}x_{1}\frac{e^{-\eta_{1}x_{1}}}{x_{1}}\frac{e^{-0x_{12}}}{x_{12}}\nonumber \\
 & = & \frac{1}{\pi}\int_{0}^{\infty}x_{1}^{2}dx_{1}e^{-\eta_{1}x_{1}}\int_{0}^{2\pi}d\varphi\int_{0}^{\pi}d\left(cos\theta\right)\int_{0}^{\infty}\frac{\zeta_{1}^{-1/2}}{\left(\left(\zeta_{1}+1\right)x_{1}^{2}-2\zeta_{1}x_{2}x_{1}cos\theta+\zeta_{1}x_{2}^{2}\right)}d\zeta_{1}\nonumber \\
 & = & 2\int_{0}^{\infty}x_{1}^{2}dx_{1}e^{-\eta_{1}x_{1}}\int_{-1}^{1}dy\int_{0}^{\infty}\frac{\zeta_{1}^{-1/2}}{\left(\left(\zeta_{1}+1\right)x_{1}^{2}-2\zeta_{1}x_{2}x_{1}y+\zeta_{1}x_{2}^{2}\right)}d\zeta_{1}\nonumber \\
 & = & 2\int_{0}^{\infty}x_{1}^{2}dx_{1}e^{-\eta_{1}x_{1}}\int_{0}^{\infty}\zeta_{1}^{-1/2}d\zeta_{1}\nonumber \\
 & \times & \frac{\log\left(\left(\zeta_{1}+1\right)x_{1}^{2}+2\zeta_{1}x_{2}x_{1}+\zeta_{1}x_{2}^{2}\right)-\log\left(\left(\zeta_{1}+1\right)x_{1}^{2}-2\zeta_{1}x_{2}x_{1}+\zeta_{1}x_{2}^{2}\right)}{2\zeta_{1}x_{1}x_{2}}\nonumber \\
 & = & 2\int da\int_{0}^{\infty}x_{1}^{2}dx_{1}e^{-\eta_{1}x_{1}}\int_{0}^{\infty}\zeta_{1}^{-1/2}d\zeta_{1}\nonumber \\
 & \times & \left(\frac{1}{2x_{1}x_{2}\left(a\zeta_{1}+\left(\zeta_{1}+1\right)x_{1}^{2}+2\zeta_{1}x_{2}x_{1}\right)}-\frac{1}{2x_{1}x_{2}\left(at_{1}+\left(\zeta_{1}+1\right)x_{1}^{2}-2\zeta_{1}x_{2}x_{1}\right)}\right)+C\nonumber \\
 & = & \frac{4\pi\left(1-e^{-\eta_{1}x_{2}}\right)}{x_{2}\eta_{1}^{2}}+C\quad,
\end{eqnarray}
 where we had to ``renormalize'' \cite{Harris Frolov and Smith JCP 120 9974 (2004)}
this infinite logarithmic integral by taking its derivative with respect
to $a=x_{2}^{2}$, whereupon integration over \emph{t} was possible
followed by \emph{a} and then $x_{1}$. In the $\zeta_{1}$ integral,
if we set $x_{2}\rightarrow\infty$ the integral goes to zero. But
this is also true in the last line above only if $C=0$, giving the
correct limit of eq. (\ref{eq:finalresult}). We see in this sequence
that a simpler integral-transform does not necessarily lead to an
easier flow. But one can hope for both. In addition, the main failing
of this integral transform is that it does not seem to be generalizable
to Slater orbitals and, hence, is of lesser value.

In an exploration of alternatives, one notes that there are also several
other integral transforms that might take the place of the Fourier
approach and involve one-dimensional integrals rather than three per
wave function. Consider for example application of \cite{GR5 p. 467 No. 3.773. GR7 p. 444}

\begin{equation}
\frac{e^{-\eta_{1}x_{1}}}{x_{1}}\frac{e^{-\eta_{12}x_{12}}}{x_{12}}=\int_{0}^{\infty}\int_{0}^{\infty}\frac{2}{\pi}\frac{\cos\left(t_{1}\eta_{1}\right)}{\left(t_{1}^{2}+x_{1}^{2}\right)}\frac{2}{\pi}\frac{\cos\left(t_{2}\eta_{12}\right)}{\left(t_{2}^{2}+x_{12}^{2}\right)}dt_{1}dt_{2}\label{eq:zt1t2_trans}
\end{equation}
 to the initial problem.

We may again use eq. (\ref{eq:GR5 p. 649 No. 4.638.2 GR7 p. 615})
(with $ $$p_{1}=1$ and $s=2$) to move both denominators into a
common quadratic form

\begin{equation}
\frac{e^{-\eta_{1}x_{1}}}{x_{1}}\frac{e^{-\eta_{12}x_{12}}}{x_{12}}=\frac{4}{\pi^{2}}\int_{0}^{\infty}\int_{0}^{\infty}\cos\left(t_{1}\eta_{1}\right)\cos\left(t_{2}\eta_{12}\right)\int_{0}^{\infty}\frac{\zeta_{1}^{1-1}}{\left(\zeta_{1}\left(t_{1}^{2}+x_{1}^{2}\right)+\left(t_{2}^{2}+x_{12}^{2}\right)\right)^{2}}d\zeta_{1}dt_{1}dt_{2}\label{eq:t1t2zeta-1}
\end{equation}

As the author was acknowledging that this was no improvement on the
Gaussian Transform, a creative flash led to the following question:
``What happens if instead of completing the square in the coordinate
variables and integrating, one does the $t$ integrals first?'' The
$t_{2}$ integral is just \cite{GR5 p. 467 No. 3.773. GR7 p. 444}
again so that, \cite{GR5 p. 772 No. 6.726.4 GR7 p .738}

\begin{align}
\frac{e^{-\eta_{1}x_{1}}}{x_{1}}\frac{e^{-\eta_{12}x_{12}}}{x_{12}} & =\frac{2}{\pi}\int_{0}^{\infty}\int_{0}^{\infty}\cos\left(t_{1}\eta_{1}\right)\frac{2}{\pi}\sqrt{\pi}2^{-3/2}\eta_{12}^{3/2}\frac{K_{\frac{3}{2}}\left(\eta_{12}\sqrt{x_{12}^{2}+\left(t_{1}^{2}+x_{1}^{2}\right)\zeta_{1}}\right)}{\left(\sqrt{x_{12}^{2}+\left(t_{1}^{2}+x_{1}^{2}\right)\zeta_{1}}\right){}^{3/2}}d\zeta_{1}dt_{1}\nonumber \\
 & =\int_{0}^{\infty}\frac{\sqrt{\zeta_{1}\eta_{12}^{2}+\eta_{1}^{2}}K_{1}\left(\sqrt{x_{1}^{2}+\frac{x_{12}^{2}}{\zeta_{1}}}\sqrt{\eta_{1}^{2}+\zeta_{1}\eta_{12}^{2}}\right)}{\pi\zeta_{1}^{3/2}\sqrt{x_{1}^{2}+\frac{x_{12}^{2}}{\zeta_{1}}}}d\zeta_{1}\quad.\label{eq:YY}
\end{align}

This result is nothing less than the desired integral transform to take the place of
eq. (\ref{eq:GR5 p. 649 No. 4.638.2 GR7 p. 615}) for the case of
two Slater orbitals. We will see that the fact that the quadratic
form $x_{1}^{2}+\frac{x_{12}^{2}}{\zeta_{1}}$ appears in two places
is not an impediment. One simply completes the square and copies the
result from one quadratic form to its identical mate. So this may
indeed be simpler than the Gaussian transform for some problems.

Having the desired integral transform in hand, let us apply it to
the original problem. %
{} We first complete the square in the quadratic form (changing variables
from $\mathbf{x}_{1}$ to $\mathbf{x}'_{1}=\mathbf{x}_{1}-\frac{1}{\zeta_{1}+1}\mathbf{x}_{2}$
with unit Jacobian) so that \cite{GR5 p. 727 No. 6.596.3 GR7 p. 693,GR5 p. 111 No. 2.311}

\begin{eqnarray}
S_{1}^{\eta_{1}0\eta_{12}0}\left(0;0,x_{2}\right) & = & \int d^{3}x_{1}\int_{0}^{\infty}\frac{\sqrt{\zeta_{1}\eta_{12}^{2}+\eta_{1}^{2}}K_{1}\left(\sqrt{x_{1}^{2}+\frac{x_{12}^{2}}{\zeta_{1}}}\sqrt{\eta_{1}^{2}+\zeta_{1}\eta_{12}^{2}}\right)}{\pi\zeta_{1}^{3/2}\sqrt{\frac{x_{12}^{2}}{\zeta_{1}}+x_{1}^{2}}}d\zeta_{1}\nonumber \\
 & = & \int d^{3}x'_{1}\int_{0}^{\infty}\frac{\sqrt{\zeta_{1}\eta_{12}^{2}+\eta_{1}^{2}}K_{1}\left(\sqrt{\frac{x_{1}^{' 2}\left(\zeta_{1}+1\right)}{\zeta_{1}}+\frac{x_{2}^{2}}{\zeta_{1}+1}}\sqrt{\eta_{1}^{2}+\zeta_{1}\eta_{12}^{2}}\right)}{\pi\zeta_{1}^{3/2}\sqrt{\frac{x_{1}^{' 2}\left(\zeta_{1}+1\right)}{\zeta_{1}}+\frac{x_{2}^{2}}{\zeta_{1}+1}}}d\zeta_{1}\nonumber \\
 & = & \int_{0}^{\infty}\frac{2\pi e^{-\frac{x_{2}\sqrt{\zeta_{1}\eta_{12}^{2}+\eta_{1}^{2}}}{\sqrt{\zeta_{1}+1}}}}{\left(\zeta_{1}+1\right){}^{3/2}\sqrt{\zeta_{1}\eta_{12}^{2}+\eta_{1}^{2}}}d\zeta_{1}=\int_{x_{2}\eta_{1}}^{x_{2}\eta_{12}}\frac{4\pi e^{-y}}{x_{2}\left(\eta_{12}^{2}-\eta_{1}^{2}\right)}\,dy\nonumber \\
 & = & \frac{4\pi\left(e^{-\eta_{12}x_{2}}-e^{-\eta_{1}x_{2}}\right)}{x_{2}\left(\eta_{1}^{2}-\eta_{12}^{2}\right)}\quad,\label{eq:syyvianew}
\end{eqnarray}
 which is indeed a much shorter path to the solution than the Gaussian
and Fourier transforms give.

\section{Generalization}

In principle, any integral transform that converts a Slater orbital
into a denominator of some power -- to be combined with an integral
transform like \cite{GR5 p. 649 No. 4.638.2 GR7 p. 615} -- could
be used for multiple products of Slater orbitals, for instance the
Stieltjes Transform \cite{ET II p. 220 No. 14.2.41} or the Bessel
function equivalent of the transform in eq. (\ref{eq:zt1t2_trans}),
\cite{GR5 p. 706 No. 6.554.4 GR7 p. 675}

\begin{equation}
\frac{e^{-\lambda r}}{r}=\int_{0}^{\infty}\frac{xJ_{0}(x\lambda)}{\left(r^{2}+x^{2}\right)^{3/2}}\,dx\quad.\label{eq:J0trans}
\end{equation}

It turns out that using the Fourier transform as this stepping stone
most easily allows one to find the general form for the equivalent
integral transform of eq. (\ref{eq:YY}) for a product of M Slater
orbitals if one takes the additional step of moving the denominator
into an exponential using \cite{GR5 p. 364 No. 3.381.4} 
\begin{equation}
(\nu-1)!D^{-\nu}=\int_{0}^{\infty}d\rho\rho^{\nu-1}e^{-\rho D}\quad:\label{eq:8}
\end{equation}

\begin{eqnarray}
\frac{e^{-R_{1}\eta_{1}}}{R_{1}}\hspace{-0.25cm} & \cdot & \hspace{-0.25cm}\frac{e^{-R_{2}\eta_{2}}}{R_{2}}\cdots\frac{e^{-R_{M}\eta_{M}}}{R_{M}}=\int d^{3}k_{1}\int d^{3}k_{2}\cdots\int\,d^{3}k_{M}\frac{1}{2\pi^{2}}\cdot\frac{e^{ik_{1}\cdot R_{1}}}{k_{1}^{2}+\eta_{1}^{2}}\cdot\frac{1}{2\pi^{2}}\cdot\frac{e^{ik_{2}\cdot R_{2}}}{k_{2}^{2}+\eta_{2}^{2}}\cdots\,\frac{1}{2\pi^{2}}\cdot\frac{e^{ik_{M}\cdot R_{M}}}{k_{M}^{2}+\eta_{M}^{2}}\nonumber \\
 & = & \int_{0}^{\infty}d\zeta_{1}\int_{0}^{\infty}d\zeta_{2}\cdots\int_{0}^{\infty}\,d\zeta_{M-1}\int d^{3}k_{1}\int d^{3}k_{2}\cdots\int\,d^{3}k_{M}\nonumber \\
 & \times & \frac{(M-1)!}{16\pi^{8}}\frac{\exp\left(ik_{1}\cdot R_{1}+ik_{2}\cdot R_{2}+\cdots+ik_{M-1}\cdot x_{M-1}+ik_{M}\cdot R_{M}\right)}{\left(\left(k_{1}^{2}+\eta_{1}^{2}\right)+\zeta_{1}\left(k_{2}^{2}+\eta_{2}^{2}\right)+\zeta_{2}\left(k_{3}^{2}+\eta_{13}^{2}\right)+\cdots+\zeta_{M-1}\left(k_{M}^{2}+\eta_{M}^{2}\right)\right)^{M}}\nonumber \\
 & = & \frac{1}{2^{M}\pi^{2M}}\int_{0}^{\infty}d\rho\int_{0}^{\infty}d\zeta_{1}\int_{0}^{\infty}d\zeta_{2}\cdots\int_{0}^{\infty}\,d\zeta_{M-1}\int d^{3}k_{1}\int d^{3}k_{2}\cdots\int\,d^{3}k_{M}\nonumber \\
 & \times & \exp\left(-\rho\left(ik_{1}\cdot R_{1}/\rho-ik_{2}\cdot R_{2}/\rho-\cdots-ik_{M-1}\cdot x_{M-1}/\rho-ik_{M}\cdot R_{M}/\rho\right)\right)\nonumber \\
 & \times & \rho^{M-1}exp\left(-\rho\left(\left(k_{1}^{2}+\eta_{1}^{2}\right)+\zeta_{1}\left(k_{2}^{2}+\eta_{2}^{2}\right)+\zeta_{2}\left(k_{3}^{2}+\eta_{13}^{2}\right)+\cdots+\zeta_{M-1}\left(k_{M}^{2}+\eta_{M}^{2}\right)\right)\right)\nonumber \\
 & = & \frac{1}{2^{M}\pi^{2M}}\int_{0}^{\infty}d\rho\int_{0}^{\infty}d\zeta_{1}\int_{0}^{\infty}d\zeta_{2}\cdots\int_{0}^{\infty}\,d\zeta_{M-1}\int d^{3}k_{1}\int d^{3}k_{2}\cdots\int\,d^{3}k_{M}\exp\left(-\rho Q\right)\quad.
\end{eqnarray}

The quadratic form may be written as \cite{Stra89c}

\begin{equation}
Q=\underline{V}^{T}\underline{W}\underline{V}\label{eq:10}
\end{equation}
 where
\begin{equation}
\underline{V}^{T}=\left(\mathbf{k}_{1},\,\mathbf{k}_{2},\cdots\mathbf{,\,k}_{M},1\right)\quad,\label{eq:11}
\end{equation}

\begin{equation}
\underline{W}=\left(\begin{array}{ccccc}
1 & 0 & \cdots & 0 & {\bf b}_{1}\\
0 & \zeta_{1} & \cdots & 0 & {\bf b}_{2}\\
 \vdots\ &  \vdots\ & \ddots &  \vdots\ & \vdots\\
0 & 0 & \cdots & \zeta_{M-1} & {\bf b}_{M}\\
{\bf b}_{1} & {\bf b}_{2} & \cdots & {\bf b}_{M} & C
\end{array}\right)\;\;,\label{eq:12}
\end{equation}

\begin{equation}
C=\eta_{1}^{2}+\zeta_{1}\eta_{2}^{2}+\zeta_{2}\eta_{3}^{2}+\cdots+\zeta_{M-1}\eta_{M}^{2}\;\;,\label{eq:13}
\end{equation}
 and

\begin{equation}
\mathbf{b_{j}}=-\frac{i}{2\rho}\mathbf{R_{j}}\quad.\label{eq:15}
\end{equation}


Now suppose one could find an orthogonal transformation that reduced
Q to diagonal form

\begin{equation}
Q'=a'_{1}k{}_{1}^{2}+a'_{'2}k_{2}^{'2}+\ldots+a'_{N+M}k_{N+M}^{'2}+c',\label{eq:16}
\end{equation}
 where, as shown by Chisholm, \cite{Chisholm} the $a$' are positive.
Then after a simple translation in $\lbrace\mathbf{k}_{1},\,\mathbf{k}_{2},\cdots,\mathbf{\,k}_{M}\rbrace$
space (with Jacobian = 1), the \emph{k} integrals could be done,\cite{GR5 p. 382 No. 3.461.2}
\begin{equation}
\int\,d^{3}k'_{1}\ldots\,d^{3}k'_{M}e^{-\rho\left(a'_{1}k{}_{1}^{'2}+a'_{2}k_{2}^{'2}+\ldots+a'_{M}k_{M}^{'2}\right)}=\left(\frac{\pi^{M}}{\rho^{M}\Lambda}\right)^{3/2},\label{eq:17}
\end{equation}
 leaving just the exponential of $-\rho c'$ to integrate over $\rho$
and the $\zeta_{i}$. But $\Lambda$ is an invariant determinant

\begin{equation}
\Lambda=\left|\begin{array}{cccc}
1 & 0 &  \cdots \ & 0\\
0 & \zeta_{1} &   \cdots & 0\\
 \vdots &  \vdots & \ddots & 0\\
0 & 0 &  \cdots & \zeta_{M-1}
\end{array}\right|=\left(1\right)\prod_{i=1}^{M-1}\zeta_{i}=\prod_{i=1}^{M}a'_{i}\quad,\label{eq:18}
\end{equation}
 so actually finding the orthogonal transformation that reduces Q
to diagonal form is unnecessary.

This orthogonal transformation also leaves

\begin{equation}
\Omega=\mathrm{det}\mathbf{W}\label{eq:35}
\end{equation}

\noindent invariant and to find its value one need only expand $\Omega$
by minors:

\begin{equation}
{\displaystyle c'\Lambda=\Omega=C\Lambda+\sum_{i=1}^{M}\sum_{j=1}^{M}\mathrm{\mathbf{b}}_{i}\cdot\mathrm{\mathbf{b}}_{j}\left(-1\right)^{i+j+1}\Lambda_{ij}=C\Lambda-b_{1}^{2}\left(1\right)\prod_{i=1}^{M-1}\zeta_{i}-\sum_{j=2}^{M}b_{j}^{2}\prod_{i\neq j-1}^{M-1}\zeta_{i}}\label{eq:40}
\end{equation}

\noindent where $\Lambda_{ij}$ is $\Lambda$ with the \emph{i}th
row and \emph{j}th column deleted, and this is diagonal in the present
case. Therefore $c'$ (of eq. (\ref{eq:16})) is given by

\begin{eqnarray}
c'=\Omega/\Lambda & = & \eta_{1}^{2}+\sum_{j=2}^{M}\zeta_{j-1}\eta_{j}^{2}-b_{1}^{2}-\sum_{j=2}^{M}b_{j}^{2}\frac{1}{\zeta_{j-1}}\nonumber \\
 & = & \eta_{1}^{2}+\sum_{j=2}^{M}\zeta_{j-1}\eta_{j}^{2}+\frac{R_{1}^{2}}{4\rho^{2}}+\sum_{j=2}^{M}\frac{R_{j}^{2}}{4\rho^{2}\zeta_{j-1}}\label{eq:39}
\end{eqnarray}

so that

\begin{equation}
\begin{array}{ccc}
\frac{e^{-R_{1}\eta_{1}}}{R_{1}}\cdot\frac{e^{-R_{2}\eta_{2}}}{R_{2}}\cdots\frac{e^{-R_{M}\eta_{M}}}{R_{M}} & = & \frac{1}{2^{M}\pi^{2M}}\int_{0}^{\infty}d\rho\int_{0}^{\infty}d\zeta_{1}\int_{0}^{\infty}d\zeta_{2}\cdots\int_{0}^{\infty}\,d\zeta_{M-1}{\displaystyle \frac{\pi^{3M/2}}{\rho^{M/2+1}\prod_{i=1}^{M-1}\zeta_{i}^{3/2}}}\\
 & \times & \hspace{-1cm}exp\left(-\rho\left(\eta_{1}^{2}+\zeta_{1}\eta_{2}^{2}+\zeta_{2}\eta_{3}^{2}+\cdots\,+\zeta_{M-1}\eta_{M}^{2}\right)\right)\\
 & \times &exp\left(-{\displaystyle \left(R_{1}^{2}+\frac{R_{2}^{2}}{\zeta_{1}}+\frac{R_{3}^{2}}{\zeta_{2}}+\cdots\,+\frac{R_{M}^{2}}{\zeta_{M-1}}\right)\frac{1}{4\rho}}\right)\quad.
\end{array}\label{eq:MRpoForm}
\end{equation}
 We perform the $\rho$ integral \cite{GR5 p. 384 No. 3.471.9} to
give the most compact, final form for this integral transform:

\begin{equation}
\begin{array}{ccc}
\frac{e^{-R_{1}\eta_{1}}}{R_{1}}\hspace{-0.3cm} & \cdot & \hspace{-0.6cm}\frac{e^{-R_{2}\eta_{2}}}{R_{2}}\cdots\frac{e^{-R_{M}\eta_{M}}}{R_{M}}=\frac{1}{2^{M}\pi^{2M}}\int_{0}^{\infty}d\zeta_{1}\int_{0}^{\infty}d\zeta_{2}\cdots\int_{0}^{\infty}\,d\zeta_{M-1}{\displaystyle \frac{\pi^{3M/2}}{\prod_{i=1}^{M-1}\zeta_{i}^{3/2}}}2^{\frac{M}{2}+1}\\
 & \times & \left(R_{1}^{2}+\frac{R_{2}^{2}}{\zeta_{1}}+\frac{R_{3}^{2}}{\zeta_{2}}+\cdots\,+\frac{R_{M}^{2}}{\zeta_{M-1}}\right)^{-M/4}\left(\eta_{1}^{2}+\zeta_{1}\eta_{2}^{2}+\zeta_{2}\eta_{3}^{2}+\cdots\,+\zeta_{M-1}\eta_{M}^{2}\right)^{M/4}\\
 & \times & K_{\frac{M}{2}}\left(\sqrt{R_{1}^{2}+\frac{R_{2}^{2}}{\zeta_{1}}+\frac{R_{3}^{2}}{\zeta_{2}}+\cdots\,+\frac{R_{M}^{2}}{\zeta_{M-1}}}\sqrt{\eta_{1}^{2}+\zeta_{1}\eta_{2}^{2}+\zeta_{2}\eta_{3}^{2}+\cdots\,+\zeta_{M-1}\eta_{M}^{2}}\right)\quad.
\end{array}\label{eq:MtransCompact}
\end{equation}

There may be some problems for which having the inverse integration
variables associated with the $\eta$'s rather than the coordinate
variables would be desirable. A simple change of variables to $\zeta_{i}=\frac{1}{\xi_{i}}$
accomplishes this:

\begin{equation}
\begin{array}{ccc}
\frac{e^{-R_{1}\eta_{1}}}{R_{1}}\cdot\frac{e^{-R_{2}\eta_{2}}}{R_{2}}\cdots\frac{e^{-R_{M}\eta_{M}}}{R_{M}} & = & \frac{1}{2^{M}\pi^{2M}}\int_{0}^{\infty}d\rho\int_{0}^{\infty}d\xi_{1}\int_{0}^{\infty}d\xi_{2}\cdots\int_{0}^{\infty}\,d\xi_{M-1}{\displaystyle \frac{\pi^{3M/2}}{\rho^{M/2+1}\prod_{i=1}^{M-1}\xi_{i}^{1/2}}}\\
 & \times & \hspace{-1cm}exp\left(-\rho\left(\eta_{1}^{2}+\frac{\eta_{2}^{2}}{\xi_{1}}+\frac{\eta_{3}^{2}}{\xi_{2}}+\frac{\eta_{4}^{2}}{\xi_{3}}+\cdots\,+\frac{\eta_{M}^{2}}{\xi_{M-1}}\right)\right)\\
 & \times & exp\left(-{\displaystyle \left(R_{1}^{2}+\xi_{1}R_{2}^{2}+\xi_{2}R_{3}^{2}+\cdots\,+\xi_{M-1}R_{M}^{2}\right)\frac{1}{4\rho}}\right)\quad.
\end{array}\label{eq:MInvRhoForm}
\end{equation}

\begin{equation}
\begin{array}{ccc}
\frac{e^{-R_{1}\eta_{1}}}{R_{1}}\hspace{-0.2cm} & \cdot & \hspace{-2.3cm}\frac{e^{-R_{2}\eta_{2}}}{R_{2}}\cdots\frac{e^{-R_{M}\eta_{M}}}{R_{M}}=\frac{1}{2^{M}\pi^{2M}}\int_{0}^{\infty}d\xi_{1}\int_{0}^{\infty}d\xi_{2}\cdots\int_{0}^{\infty}\,d\xi_{M-1}{\displaystyle \frac{\pi^{3M/2}}{\prod_{i=1}^{M-1}\xi_{i}^{1/2}}}2^{\frac{M}{2}+1}\quad,\\
 & \times & \left(R_{1}^{2}+\xi_{1}R_{2}^{2}+\xi_{2}R_{3}^{2}+\cdots\,+\xi_{M-1}R_{M}^{2}\right)^{-M/4}\left(\eta_{1}^{2}+\frac{\eta_{2}^{2}}{\xi_{1}}+\frac{\eta_{3}^{2}}{\xi_{2}}+\frac{\eta_{4}^{2}}{\xi_{3}}+\cdots\,+\frac{\eta_{M}^{2}}{\xi_{M-1}}\right)^{M/4}\\
 & \times & K_{\frac{M}{2}}\left(\sqrt{R_{1}^{2}+\xi_{1}R_{2}^{2}+\xi_{2}R_{3}^{2}+\cdots\,+\xi_{M-1}R_{M}^{2}}\sqrt{\eta_{1}^{2}+\frac{\eta_{2}^{2}}{\xi_{1}}+\frac{\eta_{3}^{2}}{\xi_{2}}+\frac{\eta_{4}^{2}}{\xi_{3}}+\cdots\,+\frac{\eta_{M}^{2}}{\xi_{M-1}}}\right)\quad.
\end{array}\label{eq:MInvCompactForm}
\end{equation}

\subsection{Inclusion of plane waves and dipole interactions}

Transition amplitudes containing plane waves may be easily included
in this integral transform, either directly in the $\rho$ version
prior to completing the square -- by utilizing orthogonal transformation
that reduces the spatial-coordinate quadratic form to diagonal form,
which never needs to actually be determined, followed by a simple
translation in $\lbrace\mathbf{x}_{1},\,\mathbf{x}_{2},\cdots,\mathbf{\,x}_{N}\rbrace$
space (with Jacobian = 1) -- or in the more compact version simply
by applying the translation in $\lbrace\mathbf{x}_{1},\,\mathbf{x}_{2},\cdots,\mathbf{\,x}_{N}\rbrace$
space to the plane wave(s) that multiply eqs. (\ref{eq:MtransCompact})
and (\ref{eq:MInvCompactForm}).

Photoionization transition amplitudes will generally contain dipole
terms $\cos\left(\theta\right)$ that may be transformed into plane
waves via a transformation like \cite{Stra90a} $\cos\theta_{12}=-x_{1}^{-1}x_{2}^{-1}\left.\frac{\partial}{\partial Q}\right|_{Q=0}e^{-Q{\bf x}_{1}\cdot{\bf x}_{2}}$,
giving an integro-differential transform, whose inclusion follows
that for other sorts of plane waves.

\subsection{Recursion}

One unusual feature of this integral transform is that one may apply
the recursion relationships of Macdonald functions to lower (or raise)
the indices. In particular, every trio of Slater orbitals may  recursively
be written as an integral of a new Slater orbital since

\begin{align}
 & \frac{\left(\zeta_{1}\eta_{2}^{2}+\zeta_{2}\eta_{3}^{2}+\eta_{1}^{2}\right){}^{3/4}K_{\frac{3}{2}}\left(\sqrt{R_{1}^{2}+\frac{R_{2}^{2}}{\zeta_{1}}+\frac{R_{3}^{2}}{\zeta_{2}}}\sqrt{\eta_{1}^{2}+\zeta_{1}\eta_{2}^{2}+\zeta_{2}\eta_{3}^{2}}\right)}{\sqrt{2}\pi^{3/2}\zeta_{1}^{3/2}\zeta_{2}^{3/2}\left(\frac{R_{2}^{2}}{\zeta_{1}}+\frac{R_{3}^{2}}{\zeta_{2}}+R_{1}^{2}\right){}^{3/4}}\nonumber \\
 & =-2\left.\frac{\partial}{\partial b}\frac{\exp\left(  -\sqrt{R_{1}^{2}+\frac{R_{2}^{2}}{\zeta_{1}}+\frac{R_{3}^{2}}{\zeta_{2}}+b}
 \sqrt{\zeta_{1}\eta_{2}^{2}+\zeta_{2}\eta_{3}^{2}+\eta_{1}^{2}}   \right)}{2\pi\zeta_{1}^{3/2}\zeta_{2}^{3/2}\sqrt{R_{1}^{2}+\frac{R_{2}^{2}}{\zeta_{1}}+\frac{R_{3}^{2}}{\zeta_{2}}+b}}\right|_{b=0}\label{eq:recursion}
\end{align}

In this way, one may craft additional forms of the transformation
that may be of use. For instance, for a product of four Slater orbitals, one may
simply apply eq. (\ref{eq:MtransCompact}) to all four orbitals simultaneous,
the first form, below, or do so only for the first three, reduce the
integrand to a new Slater orbital using eq. (\ref{eq:recursion}),
and then apply eq. (\ref{eq:MtransCompact}) a second time to the
last orbital paired with this new one, the second form, below. 

\begin{equation}
\begin{array}{ccc}
\frac{e^{-R_{1}\eta_{1}}}{R_{1}}\hspace{-0.3cm} & \cdot & \hspace{-0.6cm}\frac{e^{-R_{2}\eta_{2}}}{R_{2}}\frac{e^{-R_{3}\eta_{3}}}{R_{3}}\frac{e^{-R_{4}\eta_{4}}}{R_{4}}=\frac{1}{2^{4}\pi^{8}}\int_{0}^{\infty}d\zeta_{1}\int_{0}^{\infty}d\zeta_{2}\int_{0}^{\infty}d\zeta_{3}{\displaystyle \frac{\pi^{6}}{\zeta_{1}^{3/2}\zeta_{2}^{3/2}\zeta_{3}^{3/2}}}2^{3}\\
 & \times & \left(R_{1}^{2}+\frac{R_{2}^{2}}{\zeta_{1}}+\frac{R_{3}^{2}}{\zeta_{2}}+\frac{R_{4}^{2}}{\zeta_{3}}\right)^{-1}\left(\eta_{1}^{2}+\zeta_{1}\eta_{2}^{2}+\zeta_{2}\eta_{3}^{2}+\zeta_{3}\eta_{4}^{2}\right)^{1}\\
 & \times & K_{2}\left(\sqrt{R_{1}^{2}+\frac{R_{2}^{2}}{\zeta_{1}}+\frac{R_{3}^{2}}{\zeta_{2}}+\frac{R_{4}^{2}}{\zeta_{3}}}\sqrt{\eta_{1}^{2}+\zeta_{1}\eta_{2}^{2}+\zeta_{2}\eta_{3}^{2}+\zeta_{3}\eta_{4}^{2}}\right)\\
 & = & -2\int_{0}^{\infty}d\zeta_{1}\int_{0}^{\infty}d\zeta_{2}{\displaystyle \frac{1}{2\pi\zeta_{1}^{3/2}\zeta_{2}^{3/2}}}\frac{1}{2^{2}\pi^{4}}\int_{0}^{\infty}d\zeta_{3}{\displaystyle \frac{\pi^{3}}{\zeta_{3}^{3/2}}}2^{2}\\
 &  & \left.\frac{\partial}{\partial b}\frac{\sqrt{\zeta_{3}\left(\zeta_{1}\eta_{2}^{2}+\zeta_{2}\eta_{3}^{2}+\eta_{1}^{2}\right)+\eta_{4}^{2}}K_{1}\left(\sqrt{R_{4}^{2}+\frac{R_{1}^{2}+\frac{R_{2}^{2}}{\zeta_{1}}+\frac{R_{3}^{2}}{\zeta_{2}}+b}{\zeta_{3}}}\sqrt{\eta_{4}^{2}+\zeta_{3}\left(\eta_{1}^{2}+\zeta_{1}\eta_{2}^{2}+\zeta_{2}\eta_{3}^{2}\right)}\right)}{\pi\zeta_{3}^{3/2}\sqrt{R_{4}^{2}+\frac{R_{1}^{2}+\frac{R_{2}^{2}}{\zeta_{1}}+\frac{R_{3}^{2}}{\zeta_{2}}+b}{\zeta_{3}}}}\right|_{b=0}\\
 & = & -2\int_{0}^{\infty}d\zeta_{1}\int_{0}^{\infty}d\zeta_{2}{\displaystyle \frac{1}{2\pi\zeta_{1}^{3/2}\zeta_{2}^{3/2}}}\frac{1}{2^{2}\pi^{4}}\int_{0}^{\infty}d\xi_{1}{\displaystyle \frac{\pi^{3}}{\xi_{1}^{1/2}}}2^{2}\\
 &  & \left.\frac{\partial}{\partial b}\frac{\sqrt{\frac{1}{\xi_{1}}}\sqrt{\frac{\zeta_{1}\eta_{2}^{2}+\zeta_{2}\eta_{3}^{2}+\eta_{1}^{2}}{\xi_{1}}+\eta_{4}^{2}}K_{1}\left(\sqrt{R_{4}^{2}+\left(R_{1}^{2}+\frac{R_{2}^{2}}{\zeta_{1}}+\frac{R_{3}^{2}}{\zeta_{2}}+b\right)\xi_{1}}\sqrt{\eta_{4}^{2}+\frac{\eta_{1}^{2}+\zeta_{1}\eta_{2}^{2}+\zeta_{2}\eta_{3}^{2}}{\xi_{1}}}\right)}{\pi\sqrt{R_{4}^{2}+\left(R_{1}^{2}+\frac{R_{2}^{2}}{\zeta_{1}}+\frac{R_{3}^{2}}{\zeta_{2}}+b\right)\xi_{1}}}\right|_{b=0}
\end{array}\quad.\label{eq:MtransCompactM=00003D4}
\end{equation}

\noindent
One may instead apply eq. (\ref{eq:MInvCompactForm}) to the last
orbital and the new Slater orbital obtained from eq. (\ref{eq:recursion}),
the last form above. But redistributing $\xi_{1}$ from the first
square root in the Macdonald function to the second square root gives
a form quite similar to the first form, above, so it offers nothing
really new apart from lowing the index on the Macdonald function by
one.

\section{This transform may be used as a tool to generate a class of previously-untabled
integrals}

The compact form of this integral transform eq. (\ref{eq:MtransCompact})
involves integrals over a Macdonald function with complicated arguments
that are not tabled or found in the literature to the author's knowledge,
so the utility of the transform may well hinge on establishing their
values. This section lays out one path to that goal, comparing sequential
integration over the initially fewest number of Slater orbitals that
allow one to complete the square, with simultaneous integration over
larger numbers of Slater orbitals. The former approach will always
yield the easiest path to a solution, while comparing these two paths
will provide a suite of analytical solutions to these unusual integrals.
Although this introductory paper is perhaps not the place for a full
exploration of this set, we will sketch out the landscape of techniques
that yield solutions.

\subsection{The integral $S_{1}^{\eta_{1}0\eta_{12}0\eta_{2}0}\left(0,0;0,0,0\right)$}

Consider the next most difficult problem from eq. (\ref{eq:syyvianew}),
including a third unshifted Slater orbital and integrating over both
variables, whose solution is given easily by eq. (\ref{eq:syyvianew}),
{} and \cite{GR5 p. 358 No. 3.351.3}

\begin{eqnarray}
S_{1}^{\eta_{1}0\eta_{12}0\eta_{2}0}\left(0,0;0,0,0\right) & = & \int d^{3}x_{2}\int d^{3}x_{1}\frac{e^{-\eta_{1}x_{1}}}{x_{1}}\frac{e^{-\eta_{12}x_{12}}}{x_{12}}\frac{e^{-\eta_{2}x_{2}}}{x_{2}} \quad.\label{eq:S_Y1Y12Y2} \\
 & = & \int_{0}^{\infty}dx_{2}4\pi x_{2}^{2}\frac{4\pi\left(e^{-x_{2}\eta_{12}}-e^{-x_{2}\eta_{1}}\right)}{x_{2}\left(\eta_{1}^{2}-\eta_{12}^{2}\right)}\frac{e^{-\eta_{2}x_{2}}}{x_{2}}=\frac{16\pi^{2}}{\left(\eta_{1}+\eta_{2}\right)\left(\eta_{1}+\eta_{12}\right)\left(\eta_{2}+\eta_{12}\right)}\nonumber 
\end{eqnarray}

\subsubsection{Transforming all three Slater orbitals simultaneously}

In comparison, we next take the harder road of applying the integral
transform eq. (\ref{eq:MtransCompact}) to all three Slater orbitals
simultaneously. The integral becomes \cite{GR5 p. 727 No. 6.596.3 GR7 p. 693,GR7 p. 665 No. 6.521.10}

\begin{align}
S_{1}^{\eta_{1}0\eta_{12}0\eta_{2}0}\left(0,0;0,0,0\right) & =\int d^{3}x_{2}\int d^{3}x_{1}\frac{e^{-\eta_{1}x_{1}}}{x_{1}}\frac{e^{-\eta_{12}x_{12}}}{x_{12}}\frac{e^{-\eta_{2}x_{2}}}{x_{2}}\nonumber \\
 & =\int d^{3}x_{2}\int d^{3}x_{1}\int_{0}^{\infty}d\zeta_{1}\int_{0}^{\infty}d\zeta_{2}\nonumber \\
 & \times\frac{\left(\zeta_{1}\eta_{12}^{2}+\zeta_{2}\eta_{13}^{2}+\eta_{1}^{2}\right){}^{3/4}K_{\frac{3}{2}}\left(\sqrt{x_{1}^{2}+\frac{x_{12}^{2}}{\zeta_{1}}+\frac{x_{2}^{2}}{\zeta_{2}}}\sqrt{\eta_{1}^{2}+\zeta_{1}\eta_{12}^{2}+\zeta_{2}\eta_{2}^{2}}\right)}{\sqrt{2}\pi^{3/2}\zeta_{1}^{3/2}\zeta_{2}^{3/2}\left(\frac{x_{12}^{2}}{\zeta_{1}}+\frac{x_{2}^{2}}{\zeta_{2}}+x_{1}^{2}\right){}^{3/4}}\nonumber \\
 & =\int d^{3}x_{2}\int d^{3}x'_{1}\int_{0}^{\infty}d\zeta_{1}\int_{0}^{\infty}d\zeta_{2}\nonumber \\
 & \times\frac{\left(\zeta_{1}\eta_{12}^{2}+\zeta_{2}\eta_{13}^{2}+\eta_{1}^{2}\right){}^{3/4}K_{\frac{3}{2}}\left(\sqrt{\frac{x_{1}^{' 2}\left(\zeta_{1}+1\right)}{\zeta_{1}}+\frac{x_{2}^{2}\left(\zeta_{1}+\zeta_{2}+1\right)}{\left(\zeta_{1}+1\right)\zeta_{2}}}\sqrt{\eta_{1}^{2}+\zeta_{1}\eta_{12}^{2}+\zeta_{2}\eta_{2}^{2}}\right)}{\sqrt{2}\pi^{3/2}\zeta_{1}^{3/2}\zeta_{2}^{3/2}\left(\sqrt{\frac{x_{1}^{' 2}\left(\zeta_{1}+1\right)}{\zeta_{1}}+\frac{x_{2}^{2}\left(\zeta_{1}+\zeta_{2}+1\right)}{\left(\zeta_{1}+1\right)\zeta_{2}}}\right){}^{3/4}}\nonumber \\
 & =\int_{0}^{\infty}dx_{2}\int_{0}^{\infty}d\zeta_{1}\int_{0}^{\infty}d\zeta_{2}\frac{8\pi x_{2}^{2}}{\left(\zeta_{1}+1\right){}^{3/2}\zeta_{2}^{3/2}}K_{0}\left(\frac{x_{2}\sqrt{\zeta_{1}+\zeta_{2}+1}\sqrt{\eta_{1}^{2}+\zeta_{1}\eta_{12}^{2}+\zeta_{2}\eta_{2}^{2}}}{\sqrt{\zeta_{1}+1}\sqrt{\zeta_{2}}}\right)\nonumber \\
 & =\int_{0}^{\infty}d\zeta_{1}\int_{0}^{\infty}d\zeta_{2}\frac{4\pi^{2}}{\left(\zeta_{1}+\zeta_{2}+1\right){}^{3/2}\left(\zeta_{1}\eta_{12}^{2}+\zeta_{2}\eta_{13}^{2}+\eta_{1}^{2}\right){}^{3/2}}\nonumber \\
 & =\int_{0}^{\infty}d\zeta_{1}\left(\frac{8\pi^{2}\sqrt{\zeta_{1}+1}\eta_{13}^{2}}{\sqrt{\zeta_{1}\eta_{12}^{2}+\eta_{1}^{2}}\left(\zeta_{1}\left(\eta_{13}^{2}-\eta_{12}^{2}\right)-\eta_{1}^{2}+\eta_{13}^{2}\right){}^{2}}-\frac{16\pi^{2}\eta_{13}}{\left(\zeta_{1}\left(\eta_{13}^{2}-\eta_{12}^{2}\right)-\eta_{1}^{2}+\eta_{13}^{2}\right){}^{2}}\right.\nonumber \\
 & +\left.\frac{8\pi^{2}\sqrt{\zeta_{1}\eta_{12}^{2}+\eta_{1}^{2}}}{\sqrt{\zeta_{1}+1}\left(\zeta_{1}\left(\eta_{12}^{2}-\eta_{13}^{2}\right)+\eta_{1}^{2}-\eta_{13}^{2}\right){}^{2}}\right)\quad.\label{eq:Y1Y12Y2Simultaneous}
\end{align}
 The first and third terms of the final integral do not seem to be
tabled but the computer algebra and calculus program Mathematica 7
can do these integrals,

\begin{align}
 & \int\frac{\sqrt{a+gx}}{\sqrt{b+hx}(c+fx)^{2}}\,\,dx=\frac{1}{2}\left(\frac{2\sqrt{a+gx}\sqrt{b+hx}}{(c+fx)(ch-bf)}+\frac{(bg-ah)\log\left((c+fx)(ah-bg)\sqrt{af-cg}\sqrt{bf-ch}\right)}{\sqrt{af-cg}(bf-ch)^{3/2}}\right.\label{eq:threepolymatha7}\\
 & \left.+\frac{(ah-bg)}{\sqrt{af-cg}(bf-ch)^{3/2}}\right.\nonumber \\
 & \times\left.{\log\left(-2f(bf-ch)\left(2\sqrt{a+gx}\sqrt{b+hx}\sqrt{af-cg}\sqrt{bf-ch}+a(2bf-ch+fhx)-bcg+bfgx-2cghx\right)\right)}\right)\nonumber 
\end{align}
 yielding the result of eq. (\ref{eq:S_Y1Y12Y2}).

\subsubsection{Can one do the $\zeta_{2}$ integral before the $x_{2}$ integral?}

A more challenging question, and one quite useful to the utility of
future work, is whether one can do the $\zeta_{2}$ integral in the
fourth line of (\ref{eq:Y1Y12Y2Simultaneous}) before the $x_{2}$
integral, to generate the third line of eq. (\ref{eq:syyvianew})
were it and the new Slater orbital $\frac{e^{-\eta_{2}x_{2}}}{x_{2}}$
to be integrated over $x_{2}$. One may rewrite the Macdonald function
in terms of a Meijer G-function as \cite{PBM3 p. 665 No. 8.4.23.1}
\begin{equation}
\frac{1}{\zeta_{2}^{3/2}}K_{0}\left(2\frac{x_{2}\sqrt{\zeta_{1}+\zeta_{2}+1}\eta_{2}\sqrt{\frac{\zeta_{1}\eta_{12}^{2}+\eta_{1}^{2}}{4\eta_{2}^{2}}+\frac{\zeta_{2}}{4}}}{\sqrt{\zeta_{1}+1}\sqrt{\zeta_{2}}}\right)=\frac{1}{2}\frac{1}{\zeta_{2}^{3/2}}G_{0,2}^{2,0}\left(\frac{x_{2}^{2}\left(\zeta_{1}+\zeta_{2}+1\right)\eta_{2}^{2}\left(\zeta_{2}+\frac{\eta_{1}^{2}+\zeta_{1}\eta_{12}^{2}}{\eta_{2}^{2}}\right)}{\left(\zeta_{1}+1\right)\zeta_{2}}|\begin{array}{c}
0,0\end{array}\right)\quad.\label{eq:K0 as G2002}
\end{equation}
 One would like to use the one tabled integral \cite{PBM3 p. 349 No. 2.24.2.9}
that has roughly the right form (with $\zeta_{2}=x$),

\begin{align}
\int_{0}^{\infty}x^{\alpha-1}\left(ax^{2}+bx+c\right)^{\frac{3}{2}-\alpha} & G_{0,2}^{2,0}\left(\frac{ax^{2}+bx+c}{x}|\begin{array}{c}
\nu,-\nu\end{array}\right)\,dx\nonumber \\
 & =\frac{\sqrt{\pi}G_{1,3}^{3,0}\left(b+2\sqrt{a}\sqrt{c}|\begin{array}{c}
\frac{3}{2}\\
0,-\alpha-\nu+3,-\alpha+\nu+3
\end{array}\right)}{2a^{3/2}}\quad.\label{eq:PBM3 p. 349 No. 2.24.2.9}\\
 & +\frac{\sqrt{\pi}\sqrt{c}G_{1,3}^{3,0}\left(b+2\sqrt{a}\sqrt{c}|\begin{array}{c}
\frac{1}{2}\\
0,-\alpha-\nu+2,-\alpha+\nu+2
\end{array}\right)}{a}
\end{align}
 but inserting $\alpha=\frac{3}{2}$ to remove the polynomial multiplying $G$ in the integrand leaves us with the wrong power of \emph{x}.
One may, however, take derivatives with respect to \emph{c} of the
integrand and resultant, with $\nu=1/2$ in combination with $\nu=0$,
to show that %

\begin{align}
\int_{0}^{\infty}\frac{K_{0}\left(2\sqrt{\frac{ax^{2}+bx+c}{x}}\right)}{x^{3/2}}\,dx & =\int_{0}^{\infty}\frac{1}{x^{3/2}}\sqrt{\pi}e^{-2\sqrt{\frac{ax^{2}+bx+c}{x}}}U\left(\frac{1}{2},1,4\sqrt{\frac{ax^{2}+bx+c}{x}}\right)\,dx\nonumber \\
 & =\int_{0}^{\infty}\frac{1}{2x^{3/2}}G_{0,2}^{2,0}\left(\frac{\text{ax}^{2}+\text{b x}+c}{x}|\begin{array}{c}
0,0\end{array}\right)\,dx\nonumber \\
 & =\frac{\pi e^{-2\sqrt{2\sqrt{a}\sqrt{c}+b}}}{2\sqrt{c}\sqrt{2\sqrt{a}\sqrt{c}+b}}-\frac{\sqrt{\pi}G_{1,3}^{2,1}\left(b+2\sqrt{a}\sqrt{c}|\begin{array}{c}
-\frac{3}{2}\\
-\frac{1}{2},0,-\frac{1}{2}
\end{array}\right)}{2\sqrt{c}}\quad.\label{eq:K0/x^3/2}\\
 & =\frac{\pi e^{-2\sqrt{2\sqrt{a}\sqrt{c}+b}}}{2\sqrt{c}}\nonumber \\
 & \Rightarrow\frac{\pi e^{-2\left(\frac{2\sqrt{a}\sqrt{c}}{x_{2}\eta_{2}}+\frac{x_{2}\eta_{2}}{2}\right)}}{2\sqrt{c}}\nonumber 
\end{align}
 where the reduction of the Meijer G-function is from \cite{functions.wolfram.com/07.34.03.0727.01}
and the last step 
\begin{equation}
\sqrt{2\sqrt{a}\sqrt{c}+b}\rightarrow\frac{2\sqrt{a}\sqrt{c}}{x_{2}\eta_{2}}+\frac{x_{2}\eta_{2}}{2}\label{eq:undoing sqrt}
\end{equation}
 holds for a number of cases akin to the present one in which 
\begin{equation}
\left\{ a\to\frac{x_{2}^{2}\eta_{2}^{2}}{4\left(\zeta_{1}+1\right)},\,b\to\frac{x_{2}^{2}\eta_{2}^{2}}{4\left(\zeta_{1}+1\right)}\left(\frac{\zeta_{1}\eta_{12}^{2}+\eta_{1}^{2}}{\eta_{2}^{2}}+\zeta_{1}+1\right),\,c\to\frac{1}{4}x_{2}^{2}\left(\zeta_{1}\eta_{12}^{2}+\eta_{1}^{2}\right)\right\} \label{eq:abc->}
\end{equation}
 Inserting $\frac{2\sqrt{a}\sqrt{c}}{x_{2}\eta_{2}}+\frac{x_{2}\eta_{2}}{2}\rightarrow\frac{x_{2}\sqrt{\zeta_{1}\eta_{12}^{2}+\eta_{1}^{2}}}{2\sqrt{\zeta_{1}+1}}+\frac{x_{2}\eta_{2}}{2}$
indeed gives the third line of eq. (\ref{eq:syyvianew}) were it and
the new Slater orbital $\frac{e^{-\eta_{2}x_{2}}}{x_{2}}$ to be integrated
over $x_{2}$.

\subsubsection{Can one do the $\zeta_{2}$ integral before either the $x_{2}$ or
$x_{'1}$ integrals?}

We turn next to the question of whether the third line of eq. (\ref{eq:Y1Y12Y2Simultaneous})
can be first integrated over $\zeta_{2}$ before either of the coordinate
variables to yield the second line of eq. (\ref{eq:syyvianew}) were
it and the new Slater orbital $\frac{e^{-\eta_{2}x_{2}}}{x_{2}}$
to be integrated over $x_{2}$. %
{} The technique of the previous subsection runs into an immediate roadblock
since 
\begin{eqnarray}
\frac{\left(\zeta_{2}\eta_{2}^{2}+\zeta_{1}\eta_{12}^{2}+\eta_{1}^{2}\right){}^{3/4}}{\zeta_{2}^{3/2}\left(\frac{x_{2}^{2}\left(\zeta_{1}+\zeta_{2}+1\right)}{\left(\zeta_{1}+1\right)\zeta_{2}}+\frac{x_{1}^{' 2}\left(\zeta_{1}+1\right)}{\zeta_{1}}\right){}^{3/4}}\hspace{-1cm} &  & K_{\frac{3}{2}}\left(\sqrt{\frac{x_{1}^{' 2}\left(\zeta_{1}+1\right)}{\zeta_{1}}+\frac{x_{2}^{2}\left(\zeta_{1}+\zeta_{2}+1\right)}{\left(\zeta_{1}+1\right)\zeta_{2}}}\sqrt{\eta_{1}^{2}+\zeta_{2}\eta_{2}^{2}+\zeta_{1}\eta_{12}^{2}}\right)\nonumber \\
 & = & \frac{\left(\zeta_{2}\eta_{2}^{2}+\zeta_{1}\eta_{12}^{2}+\eta_{1}^{2}\right){}^{3/2}}{2\zeta_{2}^{3/2}\left(\left(\frac{x_{1}^{' 2}\left(\zeta_{1}+1\right)}{\zeta_{1}}+\frac{x_{2}^{2}\left(\zeta_{1}+\zeta_{2}+1\right)}{\left(\zeta_{1}+1\right)\zeta_{2}}\right)\left(\eta_{1}^{2}+\zeta_{2}\eta_{2}^{2}+\zeta_{1}\eta_{12}^{2}\right)\right){}^{3/4}})\\
 & \times & G_{0,2}^{2,0}\left(\frac{1}{4}\left(\frac{x_{1}^{' 2}\left(\zeta_{1}+1\right)}{\zeta_{1}}+\frac{x_{2}^{2}\left(\zeta_{1}+\zeta_{2}+1\right)}{\left(\zeta_{1}+1\right)\zeta_{2}}\right)\left(\eta_{1}^{2}+\zeta_{2}\eta_{2}^{2}+\zeta_{1}\eta_{12}^{2}\right)|\begin{array}{c}
\frac{3}{4},-\frac{3}{4}\end{array}\right)\nonumber 
\end{eqnarray}
 has the additional factor $\left(\zeta_{2}\eta_{2}^{2}+\zeta_{1}\eta_{12}^{2}+\eta_{1}^{2}\right){}^{3/2}$
that stands in the way of using derivatives of eq. (\ref{eq:PBM3 p. 349 No. 2.24.2.9}).
So we will utilize the $\rho$-form of the integral transform. The
$\zeta_{2}$ integral is straightforward, \cite{GR5 p. 384 No. 3.471.9}%

\begin{align}
S_{1}^{\eta_{1}0\eta_{12}0\eta_{2}0}\left(0,0;0,0,0\right) & =\int d^{3}x_{2}\int d^{3}x_{1}\frac{e^{-\eta_{1}x_{1}}}{x_{1}}\frac{e^{-\eta_{12}x_{12}}}{x_{12}}\frac{e^{-\eta_{2}x_{2}}}{x_{2}}\nonumber \\
 & =\int d^{3}x_{2}\int d^{3}x_{1}\int_{0}^{\infty}d\zeta_{1}\int_{0}^{\infty}d\zeta_{2}\int_{0}^{\infty}d\rho\nonumber \\
 & \times\frac{1}{8\pi^{3/2}\zeta_{1}^{3/2}\zeta_{2}^{3/2}\rho^{5/2}}\exp\left(-\left(\frac{x_{12}^{2}}{\zeta_{1}}+\frac{x_{2}^{2}}{\zeta_{2}}+x_{1}^{2}\right)\frac{1}{4\rho}-\rho\left(\zeta_{2}\eta_{2}^{2}+\zeta_{1}\eta_{12}^{2}+\eta_{1}^{2}\right)\right)\nonumber \\
 & =\int d^{3}x_{2}\int d^{3}x_{'1}\int_{0}^{\infty}d\zeta_{1}\int_{0}^{\infty}d\zeta_{2}\int_{0}^{\infty}d\rho\frac{1}{8\pi^{3/2}\zeta_{1}^{3/2}\zeta_{2}^{3/2}\rho^{5/2}}\nonumber \\
 & \times\exp\left(-\left(\frac{x_{1}^{' 2}\left(\zeta_{1}+1\right)}{\zeta_{1}}+\frac{x_{2}^{2}\left(\zeta_{1}+\zeta_{2}+1\right)}{\left(\zeta_{1}+1\right)\zeta_{2}}\right)\frac{1}{4\rho}-\rho\left(\zeta_{2}\eta_{2}^{2}+\zeta_{1}\eta_{12}^{2}+\eta_{1}^{2}\right)\right)\nonumber \\
 & =\int d^{3}x_{2}\int d^{3}x'_{1}\int_{0}^{\infty}d\zeta_{1}\int_{0}^{\infty}d\rho\nonumber \\
 & \times\frac{1}{4\pi x_{2}\zeta_{1}^{3/2}\rho^{2}}\exp\left(-\frac{x_{2}^{2}\zeta_{1}+\text{x1p}^{2}\left(\zeta_{1}^{2}+2\zeta_{1}+1\right)}{4\zeta_{1}\left(\zeta_{1}+1\right)\rho}-x_{2}\eta_{2}-\rho\left(\zeta_{1}\eta_{12}^{2}+\eta_{1}^{2}\right)\right)\quad.\label{eq:SYx1Yx12Yx2 zeta_2 first}
\end{align}
 as is the next integral over $\rho$ \cite{GR5 p. 384 No. 3.471.9}to
yield the second line of eq. (\ref{eq:syyvianew}) were it and the
new Slater orbital $\frac{e^{-\eta_{2}x_{2}}}{x_{2}}$ to be integrated
over $x_{2}$.

\subsection{Integrating four Slater orbitals (one shifted) over three variables, }

We will do one last test of our ability to integrate Macdonald functions
with complicated arguments by adding a fourth unshifted Slater orbital,
with the whole integrated over a third coordinate:

\begin{align}
S_{1}^{\eta_{1}0\eta_{12}0\eta_{2}0\eta_{3}0}\left(0,0,0;0,0,0,0\right) & =\int d^{3}x_{3}\int d^{3}x_{2}\int d^{3}x_{1}\frac{e^{-\eta_{1}x_{1}}}{x_{1}}\frac{e^{-\eta_{12}x_{12}}}{x_{12}}\frac{e^{-\eta_{2}x_{2}}}{x_{2}}\frac{e^{-x_{3}\eta_{3}}}{x_{3}}\nonumber \\
 & =S_{1}^{\eta_{1}0\eta_{12}0\eta_{2}0}\left(0,0;0,0,0\right)\int d^{3}x_{3}\frac{e^{-\eta_{3}x_{3}}}{x_{3}}\nonumber \\
 & =\frac{16\pi^{2}}{\left(\eta_{1}+\eta_{2}\right)\left(\eta_{1}+\eta_{12}\right)\left(\eta_{2}+\eta_{12}\right)}\int_{0}^{\infty}dx_{3}4\pi x_{3}^{2}\frac{e^{-\eta_{3}x_{3}}}{x_{3}}\nonumber \\
 & =\int_{0}^{\infty}dx_{3}4\pi x_{3}^{2}\frac{16\pi^{2}}{\left(\eta_{1}+\eta_{2}\right)\left(\eta_{1}+\eta_{12}\right)\left(\eta_{2}+\eta_{12}\right)}\frac{e^{-\eta_{3}x_{3}}}{x_{3}}\nonumber \\
 & =\frac{64\pi^{3}}{\left(\eta_{1}+\eta_{2}\right)\left(\eta_{1}+\eta_{12}\right)\left(\eta_{2}+\eta_{12}\right)\eta_{3}^{2}}\quad.\label{eq:S_Y1Y12Y2Y3}
\end{align}
 This is Harris, Frolov, and Smith's integral I (-1, -1, -1, 0, 0,
-1), \cite{Harris Frolov and Smith JCP 120 9974 (2004)} who applied
Remiddi's technique (for whom $\eta_{i}=0$) \cite{Remiddi} to simplify
the arbitrary-$\eta_{i}$ results of Fromm and Hill.\cite{Fromm and Hill}

\subsubsection{Integrating first over the \textmd{$\zeta_{3}$} variable}

We will do this in the most difficult order to generate new results:%

\begin{align}
S_{1}^{\eta_{1}0\eta_{12}0\eta_{2}0\eta_{3}0} & \left(0,0,0;0,0,0,0\right)=\int d^{3}x_{3}\int d^{3}x_{2}\int d^{3}x_{1}\frac{e^{-\eta_{1}x_{1}}}{x_{1}}\frac{e^{-\eta_{12}x_{12}}}{x_{12}}\frac{e^{-\eta_{2}x_{2}}}{x_{2}}\frac{e^{-x_{3}\eta_{3}}}{x_{3}}\nonumber \\
 & =\int d^{3}x_{3}\int d^{3}x_{2}\int d^{3}x_{1}\int_{0}^{\infty}d\zeta_{1}\int_{0}^{\infty}d\zeta_{2}\int_{0}^{\infty}d\zeta_{3}\nonumber \\
 & \times\frac{\left(\zeta_{2}\eta_{2}^{2}+\zeta_{3}\eta_{3}^{2}+\zeta_{1}\eta_{12}^{2}+\eta_{1}^{2}\right)K_{2}\left(\sqrt{x_{1}^{2}+\frac{x_{12}^{2}}{\zeta_{1}}+\frac{x_{2}^{2}}{\zeta_{2}}+\frac{x_{3}^{2}}{\zeta_{3}}}\sqrt{\eta_{1}^{2}+\zeta_{2}\eta_{2}^{2}+\zeta_{3}\eta_{3}^{2}+\zeta_{1}\eta_{12}^{2}}\right)}{2\pi^{2}\zeta_{1}^{3/2}\zeta_{2}^{3/2}\zeta_{3}^{3/2}\left(\frac{x_{12}^{2}}{\zeta_{1}}+\frac{x_{2}^{2}}{\zeta_{2}}+\frac{x_{3}^{2}}{\zeta_{3}}+x_{1}^{2}\right)}\nonumber \\
 & =\int d^{3}x_{3}\int d^{3}x_{2}\int d^{3}x'_{1}\int_{0}^{\infty}d\zeta_{1}\int_{0}^{\infty}d\zeta_{2}\int_{0}^{\infty}d\zeta_{3}\nonumber \\
 & \times\left[\frac{\zeta_{3}^{3/2}\eta_{3}^{4}}{8\pi^{2}\zeta_{1}^{3/2}\zeta_{2}^{3/2}}+\frac{\sqrt{\zeta_{3}}\eta_{3}^{2}\left(\zeta_{2}\eta_{2}^{2}+\zeta_{1}\eta_{12}^{2}+\eta_{1}^{2}\right)}{4\pi^{2}\zeta_{1}^{3/2}\zeta_{2}^{3/2}}+\frac{\left(\zeta_{2}\eta_{2}^{2}+\zeta_{1}\eta_{12}^{2}+\eta_{1}^{2}\right){}^{2}}{8\pi^{2}\zeta_{1}^{3/2}\zeta_{2}^{3/2}\sqrt{\zeta_{3}}}\right]\nonumber \\
 & \frac{K_{2}\left(2\sqrt{\frac{x_{1}^{' 2}\left(\zeta_{1}+1\right)}{4\zeta_{1}}+\frac{x_{2}^{2}\left(\zeta_{1}+\zeta_{2}+1\right)\zeta_{3}+x_{3}^{2}\left(\zeta_{1}+1\right)\zeta_{2}}{4\left(\zeta_{1}+1\right)\zeta_{2}\zeta_{3}}}\sqrt{\eta_{1}^{2}+\zeta_{2}\eta_{2}^{2}+\zeta_{3}\eta_{3}^{2}+\zeta_{1}\eta_{12}^{2}}\right)}{\left(2\left(\frac{x_{1}^{' 2}\left(\zeta_{1}+1\right)}{4\zeta_{1}}+\frac{x_{2}^{2}\left(\zeta_{1}+\zeta_{2}+1\right)\zeta_{3}+x_{3}^{2}\left(\zeta_{1}+1\right)\zeta_{2}}{4\left(\zeta_{1}+1\right)\zeta_{2}\zeta_{3}}\right)\left(\zeta_{2}\eta_{2}^{2}+\zeta_{3}\eta_{3}^{2}+\zeta_{1}\eta_{12}^{2}+\eta_{1}^{2}\right)\right)}\quad.\label{eq:SY1Y12Y2Y3}
\end{align}
 After a considerable mixing of multiple derivatives with respect
to \emph{c} and \emph{b} of the integrand and resultant of eq. (\ref{eq:PBM3 p. 349 No. 2.24.2.9})
with various values for $\nu$, we were able to determine that the
first of the three required integrals given in square brackets in
the third line above is (with $\zeta_{3}=x$) %

\begin{align}
\int_{0}^{\infty} & \frac{x^{3/2}K_{2}\left(2\sqrt{\frac{ax^{2}+bx+c}{x}}\right)}{ax^{2}+bx+c}\,dx=\int_{0}^{\infty}16\sqrt{\pi}\sqrt{x}e^{-2\sqrt{\frac{ax^{2}+bx+c}{x}}}U\left(\frac{5}{2},5,4\sqrt{\frac{ax^{2}+bx+c}{x}}\right)dx\nonumber \\
 & =\int_{0}^{\infty}\frac{x^{3/2}}{2\left(ax^{2}+bx+c\right)}G_{0,2}^{2,0}\left(\frac{ax^{2}+bx+c}{x}|\begin{array}{c}
1,-1\end{array}\right)\,dx\nonumber \\
 & =\frac{\pi e^{-2\sqrt{2\sqrt{a}\sqrt{c}+b}}}{4a^{3/2}\sqrt{2\sqrt{a}\sqrt{c}+b}}+\frac{\pi\sqrt{c}e^{-2\sqrt{2\sqrt{a}\sqrt{c}+b}}}{2a\left(2\sqrt{a}\sqrt{c}+b\right)}+\frac{\pi\sqrt{c}e^{-2\sqrt{2\sqrt{a}\sqrt{c}+b}}}{4a\left(2\sqrt{a}\sqrt{c}+b\right)^{3/2}}\nonumber \\
 & \Rightarrow e^{-x_{3}\eta_{3}}\left(\frac{\pi e^{-\frac{4\sqrt{a}\sqrt{c}}{x_{3}\eta_{3}}}}{4a^{3/2}\left(\frac{2\sqrt{a}\sqrt{c}}{x_{3}\eta_{3}}+\frac{x_{3}\eta_{3}}{2}\right)}+\frac{\pi\sqrt{c}e^{-\frac{4\sqrt{a}\sqrt{c}}{x_{3}\eta_{3}}}}{2a\left(\frac{2\sqrt{a}\sqrt{c}}{x_{3}\eta_{3}}+\frac{x_{3}\eta_{3}}{2}\right){}^{2}}+\frac{\pi\sqrt{c}e^{-\frac{4\sqrt{a}\sqrt{c}}{x_{3}\eta_{3}}}}{4a\left(\frac{2\sqrt{a}\sqrt{c}}{x_{3}\eta_{3}}+\frac{x_{3}\eta_{3}}{2}\right){}^{3}}\right)\quad.\label{eq:K2x^3/2/(ax^2+bx+c)}
\end{align}
 where the last step 
\begin{equation}
\sqrt{2\sqrt{a}\sqrt{c}+b}\rightarrow\frac{2\sqrt{a}\sqrt{c}}{x_{3}\eta_{3}}+\frac{x_{3}\eta_{3}}{2}\label{eq:undoing sqrt-1}
\end{equation}
 again holds for 
\begin{eqnarray}
\left\{ a\right. & \to & \left.\frac{1}{4}\eta_{3}^{2}\left(\frac{x_{2}^{2}\left(\zeta_{1}+\zeta_{2}+1\right)}{\left(\zeta_{1}+1\right)\zeta_{2}}+\frac{\text{x1p}^{2}\left(\zeta_{1}+1\right)}{\zeta_{1}}\right),\right.\,\nonumber \\
\left.b\right. & \to & \left.\frac{1}{4}\left(\left(\zeta_{2}\eta_{2}^{2}+\zeta_{1}\eta_{12}^{2}+\eta_{1}^{2}\right)\left(\frac{x_{2}^{2}\left(\zeta_{1}+\zeta_{2}+1\right)}{\left(\zeta_{1}+1\right)\zeta_{2}}+\frac{\text{x1p}^{2}\left(\zeta_{1}+1\right)}{\zeta_{1}}\right)+x_{3}^{2}\eta_{3}^{2}\right)\,\right.\nonumber \\
c & \to & \left.\frac{1}{4}x_{3}^{2}\left(\zeta_{2}\eta_{2}^{2}+\zeta_{1}\eta_{12}^{2}+\eta_{1}^{2}\right)\right\} \quad.\label{eq:abc->for K2}
\end{eqnarray}

Likewise, for the second term in square brackets in the third line
of (\ref{eq:SY1Y12Y2Y3}) we have

\begin{align}
\int_{0}^{\infty}\frac{\sqrt{x}}{ax^{2}+bx+c} & K_{2}\left(2\sqrt{\frac{ax^{2}+bx+c}{x}}\right)\,dx=\int_{0}^{\infty}\frac{16\sqrt{\pi}}{\sqrt{x}}e^{-2\sqrt{\frac{ax^{2}+bx+c}{x}}}U\left(\frac{5}{2},5,4\sqrt{\frac{ax^{2}+bx+c}{x}}\right)dx\nonumber \\
 & =\int_{0}^{\infty}\frac{\sqrt{x}}{2\left(ax^{2}+bx+c\right)}G_{0,2}^{2,0}\left(\frac{ax^{2}+bx+c}{x}|\begin{array}{c}
1,-1\end{array}\right)\,dx\nonumber \\
 & =\frac{\pi}{2\sqrt{a}\left(2\sqrt{a}\sqrt{c}+b\right)}e^{-2\sqrt{2\sqrt{a}\sqrt{c}+b}}\left(\frac{1}{2\sqrt{2\sqrt{a}\sqrt{c}+b}}+1\right)\nonumber \\
 & \Rightarrow\sqrt{\pi}\frac{1}{\sqrt{a}}\left(\frac{2\sqrt{a}\sqrt{c}}{x_{3}\eta_{3}}+\frac{x_{3}\eta_{3}}{2}\right){}^{-3/2}K_{\frac{3}{2}}\left(2\left(\frac{x_{3}\eta_{3}}{2}+\frac{2\sqrt{a}\sqrt{c}}{x_{3}\eta_{3}}\right)\right)\quad.\label{eq:K2x^1/2/(ax^2+bx+c)}
\end{align}

Thirdly, we have

\begin{align}
\int_{0}^{\infty}\frac{1}{\sqrt{x}\left(ax^{2}+bx+c\right)} & K_{2}\left(2\sqrt{\frac{ax^{2}+bx+c}{x}}\right)\,dx=\int_{0}^{\infty}\frac{16\sqrt{\pi}}{x^{3/2}}e^{-2\sqrt{\frac{ax^{2}+bx+c}{x}}}U\left(\frac{5}{2},5,4\sqrt{\frac{ax^{2}+bx+c}{x}}\right)dx\nonumber \\
 & =\int_{0}^{\infty}\frac{1}{2\sqrt{x}\left(ax^{2}+bx+c\right)}G_{0,2}^{2,0}\left(\frac{ax^{2}+bx+c}{x}|\begin{array}{c}
1,-1\end{array}\right)\,dx\nonumber \\
 & =-\frac{\pi e^{-2\sqrt{2\sqrt{a}\sqrt{c}+b}}}{2\sqrt{c}\left(2\sqrt{a}\sqrt{c}+b\right)}-\frac{\pi e^{-2\sqrt{2\sqrt{a}\sqrt{c}+b}}}{4\sqrt{c}\left(2\sqrt{a}\sqrt{c}+b\right)^{3/2}}\nonumber \\
 & \Rightarrow\frac{\pi e^{-2\left(\frac{2\sqrt{a}\sqrt{c}}{x_{3}\eta_{3}}+\frac{x_{3}\eta_{3}}{2}\right)}}{2\sqrt{c}\left(\frac{2\sqrt{a}\sqrt{c}}{x_{3}\eta_{3}}+\frac{x_{3}\eta_{3}}{2}\right){}^{2}}+\frac{\pi e^{-2\left(\frac{2\sqrt{a}\sqrt{c}}{x_{3}\eta_{3}}+\frac{x_{3}\eta_{3}}{2}\right)}}{4\sqrt{c}\left(\frac{2\sqrt{a}\sqrt{c}}{x_{3}\eta_{3}}+\frac{x_{3}\eta_{3}}{2}\right){}^{3}}\quad.\label{eq:K2x^-1/2/(ax^2+bx+c)}
\end{align}
 We finally sum these, weighted by their coefficients, and insert

$\frac{2\sqrt{a}\sqrt{c}}{x_{3}\eta_{3}}+\frac{x_{3}\eta_{3}}{2}\rightarrow\frac{1}{2}\sqrt{\zeta_{2}\eta_{2}^{2}+\zeta_{1}\eta_{12}^{2}+\eta_{1}^{2}}\sqrt{\frac{x_{2}^{2}\left(\zeta_{1}+\zeta_{2}+1\right)}{\left(\zeta_{1}+1\right)\zeta_{2}}+\frac{\text{x1p}^{2}\left(\zeta_{1}+1\right)}{\zeta_{1}}}+\frac{x_{3}\eta_{3}}{2}\;.$
We have checked via multidimensional numerical integration that every
step in this derivation yields the value given by the last line of
eq. (\ref{eq:S_Y1Y12Y2Y3}).

One could, of course, continue the integration process over the remaining
variables. But we will stop here since the point was not doing the
most difficult derivation of a well-known result, but exploring the
landscape of integrals over Macdonald functions with complicated arguments
that seem heretofore not to be tabled or found in the literature.

\section*{Conclusion}

We have crafted an integral transformation that may find utility in
the reduction of multidimensional transition amplitudes of quantum
theory. In particular, the general form was found for a product of
any number of Slater orbitals, whose derivatives represent hydrogenic
and Hylleraas wave functions, as well as those composed of explicitly
correlated exponentials of the kind introduced by Thakkar and Smith
\cite{Thakkar and Smith77}. In addition to atomic and nuclear transition
amplitudes, it may also find application in plasma physics, solid-state
physics, negative ion physics, and problems involving a hypothesized
non-zero-mass photon.

Unlike the Gaussian and Fourier transforms, it has the peculiarity
of displaying the quadratic form, whose square one will wish to complete,
in two locations rather than one. We have shown that this is no impediment
to its use. It has the advantage over the Gaussian transform of requiring
one fewer integrals to be subsequently reduced, and many fewer than
the $\left(3\left(M-1\right)+M-1\right)$ integral dimensions that
the Fourier transform introduces for a product of M Slater orbitals.
In cases where integrals remain, numerical integration seems to be
without problems.

Its most severe downside is likely that the quadratic forms reside
within a square root as the argument of a Macdonald function. By contrast,
Fromm and Hill \cite{Fromm and Hill} were able to leverage the nicer
form of the functions their Fourier transforms gave to integrate over
the angular and radial variables for a product of three Slater orbitals
in the three 3D integration variables with three Slater orbitals having
shifted coordinates. However, the present work is motivated by the
observation that Fromm and Hill's tour de force is unlikely to be
extensible to higher numbers of products or dimensions. Whether or
not this new approach will even approach theirs, no less exceed it,
is still an open question. We have made a start herein on finding
the analytical forms to a number of integrals over such Macdonald
functions, but as the number of functions with shifted coordinates
grows, the difficultly of doing these integrals will likely grow.
It nevertheless seems a worthwhile goal to pursue. 


\end{document}